\documentclass[11pt]{amsart}
\usepackage[T1]{fontenc}
\usepackage[utf8]{inputenc}
\usepackage{lmodern}
\usepackage{amssymb}
\usepackage[colorlinks=true,linkcolor=blue,citecolor=blue,urlcolor=blue]{hyperref}

\newtheorem{theorem}{Theorem}[section]
\newtheorem{lemma}[theorem]{Lemma}
\newtheorem{proposition}[theorem]{Proposition}
\newtheorem{corollary}[theorem]{Corollary}
\theoremstyle{definition}
\newtheorem{example}[theorem]{Example}
\theoremstyle{remark}
\newtheorem{remark}[theorem]{Remark}
\numberwithin{equation}{section}
\newcommand{\Rp}{\mathbb{R}_{+}}
\DeclareMathOperator{\Fix}{Fix}

\begin{document}

\title[Perturbed, suprametric, and interpolative contractions]{Perturbed, suprametric, and interpolative contractions are metric contractions}

\author{Firuz Kamalov}
\address{Canadian University Dubai, Dubai, United Arab Emirates}
\email{firuz@cud.ac.ae}


\begin{abstract}
We show that the fixed point theory of perturbed metric spaces, $b$-suprametric spaces, and interpolative metric spaces reduces to the metric case. A mapping is a perturbed contraction with constant $\lambda$ if and only if it is a Banach contraction with constant $\lambda$ for a metric that dominates the exact metric. Thus the theorem of Jleli and Samet is the Banach contraction principle in a finer metric, and we show that this metric cannot be chosen equivalent to the exact one. We also show that the fixed point theorem claimed for perturbed Kannan mappings is false. A $b$-suprametric is a logarithmic transform of a $b$-metric with the same constant, and an interpolative metric is a power of a metric with an explicit exponent, which is optimal for $\alpha = 1/2$. Consequently, contractive conditions of maximum type in these spaces are equivalent to the same conditions in $b$-metric and metric spaces.
\end{abstract}

\maketitle

\section{Introduction}

A large part of metric fixed point theory since the work of Banach \cite{Banach} has been concerned with relaxing the axioms of a metric space. Three recent proposals have attracted particular attention. Berzig introduced suprametric spaces \cite{Berzig1} and $b$-suprametric spaces \cite{Berzig2}, in which the triangle inequality carries an additional product term. Karap{\i}nar introduced interpolative metric spaces \cite{Kar}, in which the additional term is a weighted geometric mean. Jleli and Samet introduced perturbed metric spaces \cite{JS}, in which the observed distance $D$ is the sum of an exact metric $d$ and a nonnegative error $P$. Each of these notions has been followed by a number of papers that transfer classical contractive conditions to the new setting.

It is natural to ask how much of the resulting theory is new. For earlier generalizations of metric spaces, fixed point theorems have often turned out to be consequences of their metric counterparts after a change of metric. For the three settings above, partial results of this kind are known. Karap{\i}nar and Cvetkovi\'c \cite{KC} proved that if $d$ is a suprametric with constant $c$, then $\log(1+cd)$ is a metric. Khandani \cite{Khandani} and Karap{\i}nar and Romaguera \cite{KR} observed that interpolative metrics are $b$-metrics. Bajovi\'c, Mitrovi\'c, and Petkovi\'c \cite{BMP} showed that the continuity hypothesis in the theorem of Jleli and Samet cannot be omitted. To our knowledge, no reduction is available for perturbed metric spaces, and no transform to a $b$-metric or to a metric is known for $b$-suprametrics with $b > 1$ or for interpolative metrics.

The purpose of this paper is to give a systematic answer. We now describe the results.

Let $(X,d)$ be a metric space, let $T\colon X \to X$ be a mapping, and let $\lambda \in (0,1)$. A perturbed metric with exact metric $d$ is nothing but a function $D \ge d$; see Section~\ref{sec:prelim}. In Theorem~\ref{thm:char} we prove that $T$ satisfies $D(Tx,Ty) \le \lambda D(x,y)$ for some such $D$ if and only if
\[
\rho_\lambda(x,y) = \sup_{n \ge 0} \lambda^{-n} d(T^n x, T^n y)
\]
is finite for all $x, y \in X$, and that $\rho_\lambda$ is then a metric with $\rho_\lambda \ge d$ for which $T$ is a contraction with constant $\lambda$. This is a converse of the contraction principle in the spirit of Bessaga \cite{Bessaga}, with the additional requirement that the new metric dominate the given one. If $d$ is complete and $T$ is continuous, then $\rho_\lambda$ is complete by Theorem~\ref{thm:complete}, so the theorem of Jleli and Samet is the Banach contraction principle in $(X,\rho_\lambda)$. The metric $\rho_\lambda$ may induce a strictly finer topology than $d$, and this cannot be avoided. Indeed, Example~\ref{ex:hedgehog} exhibits a continuous perturbed contraction on a complete space that is not nonexpansive with respect to any metric topologically equivalent to $d$. Without continuity, the iterates of $T$ still converge to a point $z$ that does not depend on the initial point, and $T$ has a fixed point if and only if $Tz = z$; this is Proposition~\ref{prop:nocont}. For perturbed Kannan mappings, Example~\ref{ex:kannan} shows that Theorem~3.3 of \cite{NP} is false for every constant $\lambda \in (0,1/2)$, and Proposition~\ref{prop:kannan} gives a corrected statement.

For a $b$-suprametric $d$ with constants $b \ge 1$ and $c > 0$, we show in Theorem~\ref{thm:log} that
\[
\delta = \frac{b}{c}\log\Bigl(1 + \frac{c\,d}{b}\Bigr)
\]
is a $b$-metric with the same constant $b$. For $b = 1$ this is, up to normalization, the metric of \cite{KC}. Since $\delta$ is an increasing function of $d$, every contractive condition of the form $d(Tx,Ty) \le \theta(M(x,y))$, where $M(x,y)$ is a maximum of distances between the points $x$, $y$, $Tx$, $Ty$, is equivalent to a condition of the same form for $\delta$. This is the content of Theorem~\ref{thm:transfer}. In particular, the fixed point theorem of \cite{Berzig2} is equivalent to the corresponding theorem in $b$-metric spaces. Theorem~\ref{thm:linear} shows that a linear contraction of a $b$-suprametric space is a Banach contraction for an explicit metric.

For an $(\alpha,c)$-interpolative metric $d$, Proposition~\ref{prop:sharp} gives the $b$-metric constant $1 + c/2$, which is best possible and does not depend on $\alpha$. Theorem~\ref{thm:power} shows that $d^q$ is a metric for
\[
q = \min\Bigl\{\frac12,\ \frac{\log 2}{\log(2+c)}\Bigr\},
\]
and that this exponent cannot be improved when $\alpha = 1/2$. Thus a linear contraction of $(X,d)$ with constant $\lambda$ is a Banach contraction of $(X,d^q)$ with constant $\lambda^q$.

Finally, Section~\ref{subsec:tests} provides three tests that can be applied to any generalized triangle inequality. The paper is organized as follows. Section~\ref{sec:prelim} contains the definitions and two lemmas. Sections~\ref{sec:perturbed}, \ref{sec:supra}, and~\ref{sec:interp} treat perturbed metric spaces, $b$-suprametric spaces, and interpolative metric spaces. Section~\ref{sec:open} lists open problems, and Section~\ref{sec:conclusion} contains concluding remarks.

\section{Preliminaries}\label{sec:prelim}

Throughout the paper, $X$ is a nonempty set, $T\colon X \to X$ is a mapping, $T^0$ is the identity, and $\Fix(T)$ denotes the set of fixed points of $T$. We write $\Rp = [0,\infty)$.

\subsection{Semimetric spaces}

A \emph{semimetric} on $X$ is a function $d\colon X \times X \to \Rp$ such that $d(x,y) = 0$ if and only if $x = y$, and $d(x,y) = d(y,x)$ for all $x, y \in X$. A sequence $(x_n)$ \emph{converges} to $x$ if $d(x_n,x) \to 0$, and it is a \emph{Cauchy sequence} if $d(x_n,x_m) \to 0$ as $n, m \to \infty$. The space $(X,d)$ is \emph{complete} if every Cauchy sequence converges. We say that two semimetrics $d$ and $d'$ on $X$ are \emph{sequentially equivalent}, and we write $d \sim d'$, if they have the same convergent sequences with the same limits and the same Cauchy sequences. In this case, $(X,d)$ is complete if and only if $(X,d')$ is complete. For $\kappa \ge 0$, we say that $T$ is a \emph{$\kappa$-contraction} of $(X,d)$ if $d(Tx,Ty) \le \kappa\,d(x,y)$ for all $x, y \in X$.

Following Bessenyei and P\'ales \cite{BP}, a \emph{triangle function} for $d$ is a symmetric function $\Phi\colon [0,\infty]^2 \to [0,\infty]$ that is nondecreasing in each variable and satisfies $\Phi(0,0) = 0$ and
\begin{equation}\label{eq:triangle}
d(x,y) \le \Phi\bigl(d(x,z), d(z,y)\bigr) \qquad (x,y,z \in X).
\end{equation}
The space $(X,d)$ is \emph{regular} if it admits a triangle function that is continuous at the origin. A \emph{comparison function} is a nondecreasing function $\varphi\colon \Rp \to \Rp$ such that $\varphi^n(t) \to 0$ for every $t \ge 0$. A \emph{Matkowski function} \cite{Matkowski} is a nondecreasing function $\psi\colon \Rp \to \Rp$ such that $\psi^n(t) \to 0$ for every $t > 0$. Every Matkowski function is a comparison function. Indeed, if $\psi(t) > t$ for some $t > 0$, then $\psi^n(t) \ge t$ for all $n$ by monotonicity, which is impossible. Hence $0 \le \psi(0) \le \psi(t) \le t$ for all $t > 0$, and so $\psi(0) = 0$.

We shall use two results on regular spaces. The first is the contraction principle of \cite{BP}: if $(X,d)$ is a complete regular semimetric space and $d(Tx,Ty) \le \varphi(d(x,y))$ for all $x, y \in X$, where $\varphi$ is a comparison function, then $T$ has a unique fixed point. The second is the theorem of Chrz\k{a}szcz, Jachymski, and Turobo\'s \cite{CJT}: a semimetric space is regular if and only if it is uniformly metrizable.

\subsection{\texorpdfstring{$b$}{b}-Metrics, \texorpdfstring{$b$}{b}-suprametrics, and interpolative metrics}

A semimetric $d$ is a \emph{$b$-metric} with constant $K \ge 1$ if \eqref{eq:triangle} holds with $\Phi(u,v) = K(u+v)$. Let $p \in (0,1]$ be defined by $(2K)^p = 2$, and put
\begin{equation}\label{eq:chain}
\rho(x,y) = \inf\Bigl\{ \sum_{i=1}^{m} d(x_{i-1},x_i)^p : m \ge 1,\ x_0 = x,\ x_m = y \Bigr\}.
\end{equation}
Paluszy\'nski and Stempak \cite{PS} proved that $\rho$ is a metric on $X$ and that
\begin{equation}\label{eq:PS}
\tfrac14\, d^p \le \rho \le d^p.
\end{equation}
For $K = 1$ we have $p = 1$ and $\rho = d$.

Let $b \ge 1$ and $c \ge 0$. A semimetric $d$ is a \emph{$b$-suprametric} with constants $b$ and $c$ if \eqref{eq:triangle} holds with $\Phi(u,v) = b(u+v) + cuv$; see \cite{Berzig2}. For $b = 1$ it is a \emph{suprametric} with constant $c$; see \cite{Berzig1}.

Let $\alpha \in (0,1)$ and $c \ge 0$. Following \cite{Kar}, a semimetric $d$ is an \emph{$(\alpha,c)$-interpolative metric} if
\begin{equation}\label{eq:interp}
d(x,y) \le d(x,z) + d(z,y) + c\, d(x,z)^{\alpha} d(z,y)^{1-\alpha} \qquad (x,y,z \in X).
\end{equation}

\subsection{Perturbed metrics}

Let $D, P\colon X \times X \to \Rp$. Following Jleli and Samet \cite{JS}, the function $D$ is a \emph{perturbed metric} with respect to $P$ if $d = D - P$ is a metric on $X$. The triple $(X,D,P)$ is then called a perturbed metric space, and $d$ is called its \emph{exact metric}. Convergence, completeness, and continuity of mappings are understood with respect to $d$. We consider the conditions
\begin{equation}\label{eq:pc}
D(Tx,Ty) \le \lambda D(x,y) \qquad (x,y \in X)
\end{equation}
with $\lambda \in (0,1)$, and
\begin{equation}\label{eq:pk}
D(Tx,Ty) \le \lambda \bigl[D(x,Tx) + D(y,Ty)\bigr] \qquad (x,y \in X)
\end{equation}
with $\lambda \in [0,1/2)$. The main result of \cite{JS} states that if $(X,d)$ is complete, $T$ is continuous, and \eqref{eq:pc} holds, then $T$ has a unique fixed point. Condition \eqref{eq:pk} defines the perturbed Kannan mappings of Nu\c{t}u and P\u{a}curar \cite{NP}. The case $P = 0$ is due to Kannan \cite{Kannan}.

The definition of a perturbed metric imposes no restriction on $D$ beyond $D \ge d$. Indeed, if $d$ is a metric on $X$ and $D\colon X \times X \to \Rp$ satisfies $D \ge d$, then $P = D - d$ is nonnegative and $(X,D,P)$ is a perturbed metric space with exact metric $d$. In particular, $D$ need not be symmetric and need not vanish on the diagonal. Accordingly, from now on a perturbed metric space is a metric space $(X,d)$ together with a function $D \ge d$. We say that $T$ is a \emph{perturbed $\lambda$-contraction} of $(X,d)$ if \eqref{eq:pc} holds for some function $D\colon X \times X \to \Rp$ with $D \ge d$, and a \emph{perturbed Kannan mapping} of $(X,d)$ if \eqref{eq:pk} holds for some such $D$ and some $\lambda \in [0,1/2)$. Such a function $D$ is called \emph{admissible} for $T$.

\subsection{Conditions of maximum type}

Let $d$ be a semimetric on $X$. For $x, y \in X$ consider the five numbers
\[
d(x,y),\quad d(x,Tx),\quad d(y,Ty),\quad d(x,Ty),\quad d(y,Tx).
\]
We fix a nonempty subfamily of these numbers and denote by $M_d(x,y)$ the maximum of its members. A \emph{condition of maximum type} is a condition of the form
\[
d(Tx,Ty) \le \theta\bigl(M_d(x,y)\bigr) \qquad (x,y \in X),
\]
where $\theta\colon \Rp \to \Rp$. If it holds, we say that $T$ is a \emph{$(\theta,M)$-contraction} of $(X,d)$. The choice $M_d(x,y) = d(x,y)$ gives the condition $d(Tx,Ty) \le \theta(d(x,y))$. If $\delta$ is another semimetric on $X$, then $M_\delta$ denotes the maximum of the same subfamily formed with $\delta$ in place of $d$, and $(\theta,M)$-contractions of $(X,\delta)$ are defined accordingly.

\subsection{Two lemmas}

The first lemma converts pointwise geometric decay of the iterates into a contraction. Let $\rho$ be a metric on $X$, let $\kappa \in (0,1)$, and put
\begin{equation}\label{eq:sup}
\sigma(x,y) = \sup_{n \ge 0} \kappa^{-n} \rho\bigl(T^n x, T^n y\bigr) \qquad (x,y \in X).
\end{equation}

\begin{lemma}\label{lem:sup}
If $\sigma$ is finite, then $\sigma$ is a metric, $\sigma \ge \rho$, and $T$ is a $\kappa$-contraction of $(X,\sigma)$.
\end{lemma}

\begin{proof}
For each $n$, the function $(x,y) \mapsto \kappa^{-n}\rho(T^n x, T^n y)$ is a pseudometric on $X$. A finite pointwise supremum of pseudometrics is a pseudometric. The term with $n = 0$ shows that $\sigma \ge \rho$, and therefore $\sigma$ is a metric. Finally,
\[
\sigma(Tx,Ty) = \kappa \sup_{n \ge 1} \kappa^{-n}\rho(T^n x, T^n y) \le \kappa\,\sigma(x,y). \qedhere
\]
\end{proof}

The second lemma records that conditions of maximum type are invariant under increasing changes of the distance. Let $d$ be a semimetric on $X$, let $g\colon \Rp \to \Rp$ be an increasing bijection, and put $\delta = g \circ d$. For $\theta\colon \Rp \to \Rp$ put $\tilde\theta = g \circ \theta \circ g^{-1}$.

\begin{lemma}\label{lem:conj}
The function $\delta$ is a semimetric with $\delta \sim d$. Moreover, $T$ is a $(\theta,M)$-contraction of $(X,d)$ if and only if it is a $(\tilde\theta,M)$-contraction of $(X,\delta)$.
\end{lemma}

\begin{proof}
The functions $g$ and $g^{-1}$ are increasing and continuous, and $g(0) = 0$. Hence $\delta$ is a semimetric, $d(x_n,x) \to 0$ if and only if $\delta(x_n,x) \to 0$, and the same holds for Cauchy sequences. Since $g$ is increasing, $g(M_d(x,y)) = M_\delta(x,y)$, and the second assertion follows by applying $g$ and $g^{-1}$ to the respective inequalities.
\end{proof}

Since $\tilde\theta^{\,n} = g \circ \theta^n \circ g^{-1}$, the function $\theta$ is a Matkowski function if and only if $\tilde\theta$ is a Matkowski function.

\section{Perturbed metric spaces}\label{sec:perturbed}

Throughout this section, $(X,d)$ is a metric space, $T\colon X \to X$ is a mapping, and $\lambda \in (0,1)$. For $x, y \in X$ we put
\begin{equation}\label{eq:rho}
\rho_\lambda(x,y) = \sup_{n \ge 0} \lambda^{-n} d(T^n x, T^n y) \in [0,\infty].
\end{equation}

\subsection{Perturbed contractions are Banach contractions}

\begin{theorem}\label{thm:char}
The mapping $T$ is a perturbed $\lambda$-contraction of $(X,d)$ if and only if $\rho_\lambda$ is finite. In this case, $\rho_\lambda$ is a metric, $\rho_\lambda \ge d$, and $T$ is a $\lambda$-contraction of $(X,\rho_\lambda)$.
\end{theorem}

\begin{proof}
Let $D$ be admissible for $T$. By induction, $D(T^n x, T^n y) \le \lambda^n D(x,y)$ for all $n \ge 0$. Since $d \le D$, it follows that $\lambda^{-n} d(T^n x, T^n y) \le D(x,y)$, and hence $\rho_\lambda \le D < \infty$. Conversely, suppose that $\rho_\lambda$ is finite. By Lemma~\ref{lem:sup} with $\rho = d$ and $\kappa = \lambda$, the function $\rho_\lambda$ is a metric, $\rho_\lambda \ge d$, and $T$ is a $\lambda$-contraction of $(X,\rho_\lambda)$. In particular, $D = \rho_\lambda$ is admissible for $T$.
\end{proof}

The proof shows that $\rho_\lambda \le D$ for every admissible function $D$. Thus $\rho_\lambda$ is the smallest admissible function, and a perturbed $\lambda$-contraction of $(X,d)$ is the same as a $\lambda$-contraction with respect to some metric $\rho \ge d$.

\begin{remark}\label{rem:char}
The finiteness of $\rho_\lambda$ refers only to the metric space $(X,d)$ and the mapping $T$. It says that $d(T^n x, T^n y) = O(\lambda^n)$ for each pair of points, with a constant that may depend on the pair. The theorem of Bessaga \cite{Bessaga} provides a metric in which a given mapping is a contraction, without any relation to a given metric. Theorem~\ref{thm:char} describes the mappings for which such a metric can be found above $d$. The additive-type contractions of \cite{SSPS} are the perturbed contractions that have an admissible function $D$ such that $D - d$ is a semimetric. This requirement is no restriction, because $2\rho_\lambda$ is admissible and $2\rho_\lambda - d \ge \rho_\lambda$ is a semimetric.
\end{remark}

\begin{theorem}\label{thm:complete}
If $(X,d)$ is complete and $T$ is a continuous perturbed $\lambda$-contraction, then $(X,\rho_\lambda)$ is complete.
\end{theorem}

\begin{proof}
Let $(x_m)$ be a Cauchy sequence in $(X,\rho_\lambda)$. Since $d \le \rho_\lambda$, it is a Cauchy sequence in $(X,d)$, and so it converges in $(X,d)$ to some point $x$. Let $\varepsilon > 0$ and choose $N$ such that $\rho_\lambda(x_m,x_l) \le \varepsilon$ for all $m, l \ge N$. Fix $n \ge 0$ and $m \ge N$. Then $\lambda^{-n} d(T^n x_m, T^n x_l) \le \varepsilon$ for all $l \ge N$. Since $T^n$ is continuous, letting $l \to \infty$ we obtain $\lambda^{-n} d(T^n x_m, T^n x) \le \varepsilon$. Taking the supremum over $n$ gives $\rho_\lambda(x_m,x) \le \varepsilon$ for all $m \ge N$.
\end{proof}

By Theorems~\ref{thm:char} and~\ref{thm:complete}, the Banach contraction principle applies in $(X,\rho_\lambda)$. Hence $T$ has a unique fixed point $\xi$, and $\rho_\lambda(T^n x,\xi) \le \lambda^n(1-\lambda)^{-1}\rho_\lambda(x,Tx)$. Since $d \le \rho_\lambda \le D$ for every admissible function $D$, we obtain
\begin{equation}\label{eq:apriori}
d(T^n x, \xi) \le \frac{\lambda^n}{1-\lambda}\, \rho_\lambda(x,Tx) \le \frac{\lambda^n}{1-\lambda}\, D(x,Tx) \qquad (x \in X,\ n \ge 0).
\end{equation}

In particular, the theorem of Jleli and Samet \cite{JS} is the Banach contraction principle in the metric space $(X,\rho_\lambda)$. The right-hand side of \eqref{eq:apriori} is computable from the observed distance $D$ alone.

\subsection{The metric \texorpdfstring{$\rho_\lambda$}{rho} cannot be replaced by an equivalent metric}

By Theorem~\ref{thm:char}, we have $\rho_\lambda \ge d$, so the topology of $\rho_\lambda$ is at least as fine as that of $d$. The following example shows that it can be strictly finer, and that no other choice of metric avoids this.

\begin{example}\label{ex:hedgehog}
Let $X = \{0\} \cup \{a_{m,k} : m \ge 1,\ k \ge 0\}$, where the points are pairwise distinct. Define $w\colon X \to \Rp$ by $w(0) = 0$ and
\[
w(a_{m,k}) = 2^{k-m} \quad (0 \le k \le m), \qquad w(a_{m,k}) = \lambda^{k-m} \quad (k \ge m).
\]
Put $d(a,b) = w(a) + w(b)$ for $a \ne b$ and $d(a,a) = 0$. Define $T$ by $T0 = 0$ and $Ta_{m,k} = a_{m,k+1}$. Finally, let $N(0) = 0$ and $N(a_{m,k}) = \lambda^{k-m}$, and put $D(a,b) = N(a) + N(b)$ for $a \ne b$ and $D(a,a) = 0$. We note the following facts.
\begin{enumerate}
\item[(a)] The space $(X,d)$ is a complete metric space. Indeed, a Cauchy sequence $(y_j)$ that is not eventually constant satisfies $w(y_j) \to 0$, and hence $y_j \to 0$.
\item[(b)] The mapping $T$ is injective, and $w(Ta) \le 2w(a)$ for all $a \in X$. Hence $d(Ta,Tb) \le 2d(a,b)$, and $T$ is continuous.
\item[(c)] We have $N \ge w$ and $N(Ta) = \lambda N(a)$ for all $a \in X$. Hence $D \ge d$ and $D(Ta,Tb) = \lambda D(a,b)$ for all $a, b \in X$.
\item[(d)] Let $x_m = a_{m,0}$. Then $d(x_m,0) = 2^{-m} \to 0$, whereas $d(T^m x_m, 0) = w(a_{m,m}) = 1$ for all $m$.
\end{enumerate}
Thus the hypotheses of Theorem~\ref{thm:complete} are satisfied, and $\Fix(T) = \{0\}$. Now let $\rho'$ be a metric on $X$ that is topologically equivalent to $d$, and suppose that $\rho'(Ta,Tb) \le \rho'(a,b)$ for all $a, b \in X$. By (d), $\rho'(x_m,0) \to 0$. Hence $\rho'(T^m x_m, 0) \le \rho'(x_m,0) \to 0$, and therefore $d(T^m x_m,0) \to 0$, which contradicts (d). Consequently, $T$ is not nonexpansive, and in particular not a Banach contraction, with respect to any metric that is topologically equivalent to $d$.
\end{example}

\begin{remark}\label{rem:meyers}
Suppose that the hypotheses of Theorem~\ref{thm:complete} hold, and let $\xi$ be the fixed point of $T$. By the theorem of Meyers \cite{Meyers}, $T$ is a Banach contraction with respect to a complete metric that is topologically equivalent to $d$ if and only if the iterates $T^n$ converge to $\xi$ uniformly on some neighborhood of $\xi$. Example~\ref{ex:hedgehog} violates this condition. Thus the theorem of \cite{JS} is the Banach contraction principle in the finer metric $\rho_\lambda$, but in general not in a metric equivalent to $d$.
\end{remark}

\subsection{The role of continuity}

\begin{proposition}\label{prop:nocont}
Let $(X,d)$ be complete, and let $T$ be a perturbed $\lambda$-contraction. Then there is a point $z \in X$ such that $T^n x \to z$ for every $x \in X$, and $T$ has a fixed point if and only if $Tz = z$.
\end{proposition}

\begin{proof}
By Theorem~\ref{thm:char}, $\rho_\lambda(T^n x, T^{n+1}x) \le \lambda^n \rho_\lambda(x,Tx)$. Hence $(T^n x)$ is a Cauchy sequence in $(X,\rho_\lambda)$ and, since $d \le \rho_\lambda$, in $(X,d)$. Let $z_x$ denote its limit in $(X,d)$. For $x, y \in X$ we have $d(T^n x, T^n y) \le \lambda^n \rho_\lambda(x,y) \to 0$, and so $z_x = z_y$. Let $z$ be the common limit. If $\xi$ is a fixed point of $T$, then the sequence $(T^n \xi)$ is constant, and so $\xi = z$.
\end{proof}

Thus $z$ is the only possible fixed point of $T$. If $T$ is continuous at $z$, then $T^{n+1}x \to Tz$ and $T^{n+1}x \to z$, and so $\Fix(T) = \{z\}$.

\begin{example}[\cite{BMP}]\label{ex:BMP}
Let $X = [0,1]$ with $d(x,y) = |x-y|$, let $T0 = 1$ and $Tx = x/2$ for $x > 0$, and let $D = d + P$, where $P(x,y) = 2$ if exactly one of the points $x$, $y$ is $0$, and $P(x,y) = 0$ otherwise. Then \eqref{eq:pc} holds with $\lambda = 1/2$, and $T$ has no fixed point. In the notation of Proposition~\ref{prop:nocont}, we have $z = 0$ and $Tz = 1$. The space $(X,\rho_\lambda)$ is not complete.
\end{example}

\subsection{Perturbed Kannan mappings}

Theorem~3.3 of \cite{NP} asserts that if $(X,d)$ is complete and \eqref{eq:pk} holds, then $T$ has a unique fixed point. No continuity hypothesis is made there. The following example shows that this assertion is false for every constant $\lambda \in (0,1/2)$.

\begin{example}\label{ex:kannan}
Let $\lambda \in (0,1/2)$, and let $X = \{0\} \cup \{\lambda^n : n \ge 0\}$ with $d(x,y) = |x-y|$. The set $X$ is closed in $\mathbb{R}$, and so $(X,d)$ is complete. Define $T0 = 1$ and $T\lambda^n = \lambda^{n+1}$ for $n \ge 0$. Define a symmetric function $D\colon X \times X \to \Rp$ by
\[
D(\lambda^n,\lambda^{n+1}) = \lambda^n \quad (n \ge 0), \qquad D(0,1) = \lambda^{-1},
\]
and $D(x,y) = |x-y|$ for all other pairs. Then $D \ge d$, because $\lambda^n \ge \lambda^n - \lambda^{n+1}$ and $\lambda^{-1} \ge 1$. Moreover, $D(\lambda^n,T\lambda^n) = \lambda^n$ and $D(0,T0) = \lambda^{-1}$. We verify \eqref{eq:pk}. Let $x \ne y$. If $x = \lambda^n$ and $y = \lambda^m$ with $n < m$, then
\[
D(Tx,Ty) = D(\lambda^{n+1},\lambda^{m+1}) \le \lambda^{n+1} \le \lambda\bigl(\lambda^n + \lambda^m\bigr).
\]
If $x = 0$ and $y = \lambda^m$, then
\[
D(Tx,Ty) = D(1,\lambda^{m+1}) \le 1 = \lambda D(0,T0) \le \lambda\bigl[D(0,T0) + D(y,Ty)\bigr].
\]
For $x = y$, the left-hand side of \eqref{eq:pk} is $D(Tx,Tx) = 0$. Hence \eqref{eq:pk} holds on the complete perturbed metric space $(X,D,D-d)$, and $T$ has no fixed point.
\end{example}

The argument in \cite{NP} lets $n \to \infty$ in \eqref{eq:pk} with $x = x_n$ and $y = \lim x_n$. This step needs a continuity property of $D$, and the definition provides none. In Example~\ref{ex:kannan}, we have $T^n x \to 0$ for every $x$, but $D(T^n x,1) \to 1 \ne D(0,1)$. The correct statement is as follows.

\begin{proposition}\label{prop:kannan}
Let $(X,d)$ be complete, and let $T$ be a perturbed Kannan mapping. Then the conclusion of Proposition~\ref{prop:nocont} holds.
\end{proposition}

\begin{proof}
Let $D$ be admissible for $T$, with constant $\lambda \in [0,1/2)$, and put $\gamma = \lambda/(1-\lambda)$, so that $\gamma < 1$. Taking $y = Tx$ in \eqref{eq:pk}, we obtain the inequality $D(Tx,T^2x) \le \gamma D(x,Tx)$, and hence
\[
d(T^n x, T^{n+1}x) \le D(T^n x, T^{n+1}x) \le \gamma^n D(x,Tx).
\]
Thus $(T^n x)$ is a Cauchy sequence in $(X,d)$. Let $z_x$ be its limit. For $x, y \in X$, condition \eqref{eq:pk} gives
\[
d(T^{n+1}x, T^{n+1}y) \le D(T^{n+1}x, T^{n+1}y) \le \lambda\gamma^n \bigl[D(x,Tx) + D(y,Ty)\bigr] \to 0,
\]
and so $z_x = z_y$. The rest follows as in the proof of Proposition~\ref{prop:nocont}.
\end{proof}

As before, $\Fix(T) = \{z\}$ if $T$ is continuous at $z$.

\begin{remark}\label{rem:HR}
In Propositions~\ref{prop:nocont} and~\ref{prop:kannan}, continuity of $T$ at $z$ can be replaced by lower semicontinuity at $z$ of the function $x \mapsto d(x,Tx)$, because $d(T^n x, T^{n+1}x) \to 0$. This is the condition of Hicks and Rhoades \cite{HR}. The same proof applies to the condition $D(Tx,Ty) \le \varphi(D(x,y))$, where $\varphi$ is nondecreasing and $\sum_n \varphi^n(t) < \infty$ for every $t$. Theorem~2.2 of \cite{NP} asserts the existence of a fixed point under this condition without continuity. Example~\ref{ex:BMP} with $\varphi(t) = t/2$ shows that this is not correct, as was observed in \cite{BMP}.
\end{remark}

\section{\texorpdfstring{$b$}{b}-Suprametric spaces}\label{sec:supra}

Throughout this section, $b \ge 1$, $c > 0$, and $d$ is a $b$-suprametric on $X$ with constants $b$ and $c$. The case $c = 0$ is that of a $b$-metric. We put
\begin{equation}\label{eq:g}
g(t) = \frac{b}{c}\log\Bigl(1 + \frac{c\,t}{b}\Bigr), \qquad g^{-1}(s) = \frac{b}{c}\bigl(e^{cs/b} - 1\bigr), \qquad \delta = g \circ d.
\end{equation}

\subsection{The logarithmic transform}

\begin{theorem}\label{thm:log}
The function $\delta$ is a $b$-metric with constant $b$, and $\delta \sim d$. In particular, $\delta$ is a metric if $b = 1$.
\end{theorem}

\begin{proof}
Let $x, y, z \in X$, and put $u = d(x,z)$ and $v = d(z,y)$. Since $b \ge 1$,
\[
1 + \frac{c}{b}\, d(x,y) \le 1 + c(u+v) + \frac{c^2}{b}\,uv \le (1+cu)(1+cv) \le \Bigl(1+\frac{cu}{b}\Bigr)^{b}\Bigl(1+\frac{cv}{b}\Bigr)^{b},
\]
where the last step is the Bernoulli inequality $1 + bt \le (1+t)^b$ for $t \ge 0$. Taking logarithms and multiplying by $b/c$, we obtain $\delta(x,y) \le b[\delta(x,z) + \delta(z,y)]$. The relation $\delta \sim d$ follows from Lemma~\ref{lem:conj}.
\end{proof}

Note that $\delta \le d$, because $\log(1+t) \le t$. For $b = 1$ the theorem has a converse. If $\delta'$ is a metric on $X$ and $d' = (e^{c\delta'} - 1)/c$, then the triangle inequality for $\delta'$ gives
\[
1 + c\,d'(x,y) \le \bigl(1 + c\,d'(x,z)\bigr)\bigl(1 + c\,d'(z,y)\bigr),
\]
which is the inequality $d'(x,y) \le d'(x,z) + d'(z,y) + c\,d'(x,z)\,d'(z,y)$.

Thus the suprametrics with constant $c$ on $X$ are exactly the functions $(e^{c\delta} - 1)/c$, where $\delta$ is a metric on $X$. In other words, the suprametric inequality is the triangle inequality for $\log(1+cd)$. For $b = 1$, the metric $\delta$ is due to Karap{\i}nar and Cvetkovi\'c \cite{KC}, up to the factor $1/c$. Theorem~\ref{thm:log} extends their result to $b > 1$.

Let $p \in (0,1]$ be defined by $(2b)^p = 2$, and let $\rho$ be the function \eqref{eq:chain} formed with $\delta$ in place of $d$.

\begin{corollary}\label{cor:metrization}
The function $\rho$ is a metric on $X$, and
\begin{equation}\label{eq:rho-delta}
\tfrac14\,\delta^p \le \rho \le \delta^p.
\end{equation}
\end{corollary}

\begin{proof}
This follows from Theorem~\ref{thm:log} and \eqref{eq:PS}.
\end{proof}

For $b = 1$ we have $\rho = \delta$. By \eqref{eq:rho-delta} and Theorem~\ref{thm:log}, $\rho \sim d$. Hence the spaces $(X,d)$, $(X,\delta)$, and $(X,\rho)$ are either all complete or all incomplete.

\subsection{Contractive conditions}

Let $\theta\colon \Rp \to \Rp$, and put $\tilde\theta = g \circ \theta \circ g^{-1}$ with $g$ from \eqref{eq:g}.

\begin{theorem}\label{thm:transfer}
The mapping $T$ is a $(\theta,M)$-contraction of $(X,d)$ if and only if it is a $(\tilde\theta,M)$-contraction of $(X,\delta)$.
\end{theorem}

\begin{proof}
This is Lemma~\ref{lem:conj}.
\end{proof}

Moreover, $\theta$ is a Matkowski function if and only if $\tilde\theta$ is a Matkowski function.

It follows that every fixed point theorem of maximum type in $b$-suprametric spaces with constants $b$ and $c$ is equivalent to the same theorem in $b$-metric spaces with constant $b$. For $b = 1$, the latter are metric spaces. We record the most important case. The theorem of Berzig \cite{Berzig2} states that if $(X,d)$ is complete and $d(Tx,Ty) \le \psi(d(x,y))$ for all $x, y \in X$, where $\psi$ is a Matkowski function, then $T$ has a unique fixed point $\xi$, and $T^n x \to \xi$ for every $x \in X$.

\begin{corollary}\label{cor:berzig}
The theorem of Berzig is equivalent to the same theorem for $b$-metric spaces with constant $b$.
\end{corollary}

\begin{proof}
The equivalence follows from Theorem~\ref{thm:transfer} with $M_d(x,y) = d(x,y)$ and from Corollary~\ref{cor:metrization}.
\end{proof}

The theorem for $b$-metric spaces goes back to Czerwik \cite{Czerwik} and is contained in \cite{BP}. Indeed, a $b$-metric space is regular, since $\Phi(u,v) = b(u+v)$ is continuous at the origin, and the convergence $T^n x \to \xi$ follows from $\delta(T^n x,\xi) \le \tilde\psi^{\,n}(\delta(x,\xi))$. For $b = 1$, it is the theorem of Matkowski \cite{Matkowski} in the metric space $(X,\delta)$.

For $b = 1$, reductions of this kind are given in \cite{KC}.

\subsection{Linear contractions}

For $\theta(t) = \lambda t$ with $\lambda \in (0,1)$, the conjugate function is not linear. Indeed,
\[
\tilde\theta(s) = \frac{b}{c}\log\bigl(1 + \lambda(e^{cs/b} - 1)\bigr),
\]
so that $\tilde\theta(s)/s \to \lambda$ as $s \to 0$, whereas $s - \tilde\theta(s) \to (b/c)\log(1/\lambda)$ as $s \to \infty$. A Banach contraction is nevertheless obtained in a different metric. Let $T$ be a $\lambda$-contraction of $(X,d)$, let $p$ and $\rho$ be as in Corollary~\ref{cor:metrization}, and put
\begin{equation}\label{eq:sigma}
\sigma(x,y) = \sup_{n \ge 0} \lambda^{-np} \rho\bigl(T^n x, T^n y\bigr) \qquad (x,y \in X).
\end{equation}

\begin{theorem}\label{thm:linear}
The function $\sigma$ is a metric on $X$ with $\delta^p/4 \le \sigma \le d^p$, and $T$ is a $\lambda^p$-contraction of $(X,\sigma)$.
\end{theorem}

\begin{proof}
By \eqref{eq:rho-delta}, the inequality $\delta \le d$, and the contraction condition,
\[
\rho(T^n x, T^n y) \le \delta(T^n x, T^n y)^p \le d(T^n x, T^n y)^p \le \lambda^{np} d(x,y)^p.
\]
Hence $\sigma \le d^p < \infty$. Lemma~\ref{lem:sup} with $\kappa = \lambda^p$ shows that $\sigma$ is a metric, that $T$ is a $\lambda^p$-contraction of $(X,\sigma)$, and that $\sigma \ge \rho \ge \delta^p/4$.
\end{proof}

Since $\delta \sim d$, the theorem gives $\sigma \sim d$. In particular, $(X,\sigma)$ is complete if and only if $(X,d)$ is complete. In this case, the Banach contraction principle in $(X,\sigma)$ yields a unique fixed point $\xi$ of $T$ and
\[
\tfrac14\,\delta(T^n x,\xi)^p \le \sigma(T^n x,\xi) \le \frac{\lambda^{np}}{1-\lambda^p}\,\sigma(x,Tx) \le \frac{\lambda^{np}}{1-\lambda^p}\, d(x,Tx)^p,
\]
that is,
\begin{equation}\label{eq:apriori2}
\delta(T^n x,\xi) \le \Bigl(\frac{4}{1-\lambda^p}\Bigr)^{1/p} \lambda^n\, d(x,Tx) \qquad (x \in X,\ n \ge 0).
\end{equation}
For $b = 1$ we have $\rho = \delta$, hence $\delta \le \sigma \le d$, and the factor $4$ in \eqref{eq:apriori2} can be replaced by $1$.

Theorem~\ref{thm:linear} holds for every $\lambda \in (0,1)$, and no relation between $\lambda$, $b$, and $c$ is required. The estimate \eqref{eq:apriori2} is an explicit a priori bound, since $d(T^n x,\xi) = g^{-1}(\delta(T^n x,\xi))$.

\section{Interpolative metric spaces and general triangle functions}\label{sec:interp}

In the first two subsections, the numbers $\alpha \in (0,1)$ and $c > 0$ are fixed, and $d$ is an $(\alpha,c)$-interpolative metric on $X$. For $x, y, z \in X$ we write $u = d(x,z)$ and $v = d(z,y)$.

\subsection{The sharp \texorpdfstring{$b$}{b}-metric constant}

Khandani \cite{Khandani} proved that $d$ is a $b$-metric with constant $1 + c\max\{\alpha,1-\alpha\}$, and Karap{\i}nar and Romaguera \cite{KR} located interpolative metrics between strong $b$-metrics and $b$-metrics. The symmetry of $d$ yields a constant that does not depend on $\alpha$.

\begin{proposition}\label{prop:sharp}
For all $x, y, z \in X$,
\begin{equation}\label{eq:sym}
d(x,y) \le u + v + c\sqrt{uv} \le \Bigl(1 + \frac{c}{2}\Bigr)(u+v).
\end{equation}
The constant $1 + c/2$ is best possible for every $\alpha$ and $c$.
\end{proposition}

\begin{proof}
We apply \eqref{eq:interp} to the triples $(x,y,z)$ and $(y,x,z)$. Since $d$ is symmetric, we obtain
\[
d(x,y) \le u + v + c\min\{u^{\alpha}v^{1-\alpha},\, u^{1-\alpha}v^{\alpha}\} \le u + v + c\sqrt{uv},
\]
because the minimum of two nonnegative numbers does not exceed their geometric mean. The second inequality in \eqref{eq:sym} follows from $\sqrt{uv} \le (u+v)/2$. To see that the constant is best possible, let $X = \{x,y,z\}$ with $d(x,z) = d(z,y) = 1$ and $d(x,y) = 2 + c$. Then \eqref{eq:interp} holds, with equality for the triple $(x,y,z)$, and every $b$-metric constant $K$ of this space satisfies $2 + c \le 2K$.
\end{proof}

Thus $d$ is a $(1/2,c)$-interpolative metric and a $b$-metric with constant $1 + c/2$.

\subsection{The power transform}

\begin{lemma}\label{lem:ineq}
Let $r \ge 2$ and $s, t \ge 0$. Then
\begin{equation}\label{eq:ineq}
(s+t)^r \ge s^r + t^r + (2^r - 2)(st)^{r/2}.
\end{equation}
\end{lemma}

\begin{proof}
If $st = 0$, then \eqref{eq:ineq} is an equality. Both sides of \eqref{eq:ineq} are symmetric and homogeneous of degree $r$, so we may assume that $s = e^{\eta}$ and $t = e^{-\eta}$ with $\eta \ge 0$. Put
\[
G(\eta) = (e^{\eta} + e^{-\eta})^r - e^{r\eta} - e^{-r\eta}.
\]
Then \eqref{eq:ineq} is the inequality $G(\eta) \ge G(0)$. Let $m = r - 1 \ge 1$ and $x = e^{2\eta} \ge 1$. A direct computation gives
\[
G'(\eta) = r e^{-r\eta} \varphi_m(x), \qquad \varphi_m(x) = (x+1)^m(x-1) - x^{m+1} + 1.
\]
We show that $\varphi_m(x) \ge 0$. Clearly $\varphi_1(x) = 0$. Since $(x+1)^m \ge (x+1)x^{m-1}$, we have
\[
\frac{\partial \varphi_m}{\partial m}(x) = (x+1)^m(x-1)\log(x+1) - x^{m+1}\log x \ge x^{m-1} k(x),
\]
where $k(x) = (x^2-1)\log(x+1) - x^2\log x$. Now $k(1) = 0$ and $k'(x) = 2x\log(1 + 1/x) - 1$. The inequality $\log(1+y) \ge 2y/(2+y)$ for $y \ge 0$ gives $k'(x) \ge 4x/(2x+1) - 1 \ge 0$ for $x \ge 1$. Hence $k \ge 0$ on $[1,\infty)$, and therefore $\varphi_m \ge \varphi_1 = 0$ for all $m \ge 1$. It follows that $G' \ge 0$ on $[0,\infty)$.
\end{proof}

Let
\begin{equation}\label{eq:q}
q = \min\Bigl\{\frac12,\ \frac{\log 2}{\log(2+c)}\Bigr\}.
\end{equation}

\begin{theorem}\label{thm:power}
The function $d^q$ is a metric on $X$, and $q$ is the largest exponent with this property for all $(1/2,c)$-interpolative metrics.
\end{theorem}

\begin{proof}
By \eqref{eq:sym}, $d(x,y) \le u + v + c\sqrt{uv}$. Suppose first that $c \le 2$. Then $q = 1/2$ and $d(x,y) \le (\sqrt{u} + \sqrt{v})^2$, and so $d(x,y)^q \le u^q + v^q$. Suppose now that $c > 2$. Then $r = 1/q > 2$ and $2^r = 2 + c$, and Lemma~\ref{lem:ineq} with $s = u^q$ and $t = v^q$ gives
\[
(u^q + v^q)^{1/q} \ge u + v + c\sqrt{uv} \ge d(x,y).
\]
Thus $d^q$ satisfies the triangle inequality in both cases.

To prove the last assertion, let $0 < \varepsilon \le 1$, and let $X = \{x,y,z\}$ with $d(x,z) = \varepsilon$, $d(z,y) = 1$, and $d(x,y) = 1 + \varepsilon + c\sqrt{\varepsilon}$. This is a $(1/2,c)$-interpolative metric. Suppose that $d^s$ is a metric for all these spaces. Then
\[
(1 + \varepsilon + c\sqrt{\varepsilon})^s \le 1 + \varepsilon^s.
\]
For $\varepsilon = 1$ this gives $(2+c)^s \le 2$. Moreover, $((1 + c\sqrt{\varepsilon})^s - 1)/\sqrt{\varepsilon} \le \varepsilon^{s-1/2}$, and the left-hand side tends to $sc > 0$ as $\varepsilon \to 0$. Hence $s \le 1/2$.
\end{proof}

\begin{corollary}\label{cor:interp}
A mapping $T$ is a $\lambda$-contraction of $(X,d)$ if and only if it is a $\lambda^q$-contraction of the metric space $(X,d^q)$.
\end{corollary}

\begin{proof}
Raise the inequality $d(Tx,Ty) \le \lambda d(x,y)$ to the power $q$.
\end{proof}

More generally, Lemma~\ref{lem:conj} with $g(t) = t^q$ shows that $d^q \sim d$, and that $T$ is a $(\theta,M)$-contraction of $(X,d)$ if and only if it is a $(\tilde\theta,M)$-contraction of $(X,d^q)$, where $\tilde\theta(s) = \theta(s^{1/q})^q$. If $(X,d)$ is complete and $\lambda \in (0,1)$, then the Banach contraction principle in $(X,d^q)$ yields a unique fixed point $\xi$ of $T$ and, after raising the usual estimate to the power $1/q$,
\begin{equation}\label{eq:apriori3}
d(T^n x,\xi) \le \frac{\lambda^n}{(1-\lambda^q)^{1/q}}\, d(x,Tx) \qquad (x \in X,\ n \ge 0).
\end{equation}

Consequently, every fixed point theorem of maximum type in interpolative metric spaces is equivalent to the same theorem in the metric space $(X,d^q)$, and a linear contraction of $(X,d)$ is a Banach contraction of $(X,d^q)$.

\subsection{Three tests for a triangle function}\label{subsec:tests}

The arguments of Sections~\ref{sec:supra} and~\ref{sec:interp} used the respective axioms only through properties of their triangle functions. This leads to three tests, which apply to any axiom of the form \eqref{eq:triangle}. Let $d$ be a semimetric on $X$ with a triangle function $\Phi$. The first test is the continuity of $\Phi$ at the origin. In this case $(X,d)$ is regular, hence uniformly metrizable by \cite{CJT}, and the contraction principle of \cite{BP} applies. The second and the third test are contained in the following theorem, in which $r > 0$, $L \ge 1$, and $K \ge 1$ are constants and $g\colon \Rp \to \Rp$ is an increasing bijection.

\begin{theorem}\label{thm:tests}
\begin{enumerate}
\item[(i)] If $\Phi(t,t) \le Lt$ for all $t \in [0,r]$, then $\min\{d,r\}$ is a $b$-metric with constant $L$.
\item[(ii)] If $g(\Phi(u,v)) \le K[g(u) + g(v)]$ for all $u, v \ge 0$, then $g \circ d$ is a $b$-metric with constant $K$.
\end{enumerate}
\end{theorem}

\begin{proof}
Let $u = d(x,z)$, $v = d(z,y)$, and $d_r = \min\{d,r\}$. If $u \le r$ and $v \le r$, then the monotonicity of $\Phi$ gives
\[
d(x,y) \le \Phi(u,v) \le \Phi(\max\{u,v\},\max\{u,v\}) \le L\max\{u,v\} \le L(u+v).
\]
If $u > r$ or $v > r$, then $d_r(x,y) \le r \le L[d_r(x,z) + d_r(z,y)]$. This proves (i). Statement (ii) follows from $g(d(x,y)) \le g(\Phi(u,v)) \le K[g(u) + g(v)]$.
\end{proof}

In the situation of (i), we have $d_r \sim d$. Let $p$ be defined by $(2L)^p = 2$, and let $\rho$ be the function \eqref{eq:chain} formed with $d_r$. If $T$ is a $\lambda$-contraction of $(X,d)$ with $\lambda \in (0,1)$, then $\rho(T^n x, T^n y) \le d(T^n x, T^n y)^p \le \lambda^{np} d(x,y)^p$ by \eqref{eq:PS}. As in the proof of Theorem~\ref{thm:linear}, the function $\sigma$ defined by \eqref{eq:sigma} is then a metric with $d_r^{\,p}/4 \le \sigma \le d^p$, and $T$ is a $\lambda^p$-contraction of $(X,\sigma)$. In the situation of (ii), Lemma~\ref{lem:conj} transfers every condition of maximum type from $d$ to $g \circ d$.

Theorem~\ref{thm:log} is an instance of (ii), with $g$ from \eqref{eq:g} and $K = b$. Theorem~\ref{thm:power} is an instance of (ii), with $g(t) = t^q$, $K = 1$, and $\Phi(u,v) = u + v + c\sqrt{uv}$. Statement (i) applies to $b$-suprametrics with $L = 2b + cr$.

\begin{example}\label{ex:nonlinear}
The first test can apply when Theorem~\ref{thm:tests}(i) does not. Let $c > 0$ and $0 < \alpha < 1/2$, and let $\Phi(u,v) = u + v + c(uv)^{\alpha}$. Then $\Phi$ is continuous at the origin, and $\Phi(t,t)/t \to \infty$ as $t \to 0$. Let $X$ consist of pairwise distinct points $x_n$, $y_n$, $z_n$, where $n \ge 1$. Put $t_n = 2^{-n}$ and
\[
d(x_n,z_n) = d(z_n,y_n) = t_n, \qquad d(x_n,y_n) = 2t_n + c\,t_n^{2\alpha},
\]
and let $d(a,a') = 3 + c$ whenever $a$ and $a'$ carry different indices. Then $\Phi$ is a triangle function for $d$. Indeed, if $x$, $y$, $z$ do not carry the same index, then at least one argument of $\Phi$ in \eqref{eq:triangle} equals $3 + c$, which is the largest value of $d$. If they carry the same index, then the only nontrivial case is the triple $(x_n,y_n,z_n)$, for which \eqref{eq:triangle} is an equality. Now
\[
\frac{d(x_n,y_n)}{d(x_n,z_n) + d(z_n,y_n)} = 1 + \frac{c}{2}\,t_n^{2\alpha-1} \to \infty.
\]
Hence $\min\{d,r\}$ is not a $b$-metric for any $r > 0$, although the first test still applies.
\end{example}

\section{Open problems}\label{sec:open}

We collect some questions that the methods of this paper do not answer.

\begin{enumerate}
\item Let $q(\alpha,c)$ denote the largest $s$ such that $d^s$ is a metric for every $(\alpha,c)$-interpolative metric $d$. Theorem~\ref{thm:power} shows that $q(1/2,c) = \min\{1/2, \log 2/\log(2+c)\}$ and that $q(\alpha,c) \ge q(1/2,c)$. Three-point spaces as in the proof of Theorem~\ref{thm:power} show that
\[
q(\alpha,c) \le \min\Bigl\{\max\{\alpha,1-\alpha\},\ \frac{\log 2}{\log(2+c)}\Bigr\}.
\]
Is this upper bound attained? A random search over $2 \times 10^5$ parameter values produced no counterexample.
\item Suppose that the hypotheses of Theorem~\ref{thm:complete} hold and that $T^n \to \xi$ uniformly on a neighborhood of the fixed point $\xi$. Is $\rho_\mu$ topologically equivalent to $d$ for some $\mu \in [\lambda,1)$?
\item Let $b > 1$, and let $T$ be a Matkowski contraction of a complete $b$-metric space. Is $T$ a Matkowski contraction with respect to a sequentially equivalent metric? Theorem~\ref{thm:linear} settles the linear case.
\item Lemma~\ref{lem:conj} does not preserve conditions that involve sums, such as the conditions of Kannan, Chatterjea, and Hardy and Rogers. Are the corresponding theorems in suprametric and interpolative metric spaces consequences of metric theorems?
\item Which of the fixed point theorems for \'Ciri\'c-type and Hardy--Rogers-type conditions in perturbed metric spaces remain valid without continuity or domination hypotheses? Examples~\ref{ex:BMP} and~\ref{ex:kannan} show that the passage to the limit is the critical step.
\end{enumerate}

\section{Conclusion}\label{sec:conclusion}

The fixed point theory of the three classes of spaces considered here reduces to metric fixed point theory, in a precise sense for each class. A perturbed contraction is a Banach contraction for the metric $\rho_\lambda \ge d$, and the theorem of Jleli and Samet is the Banach contraction principle in $(X,\rho_\lambda)$. This reduction is sharp in two respects: the metric $\rho_\lambda$ cannot be replaced by a metric equivalent to $d$, and the continuity of $T$ cannot be omitted. For perturbed Kannan mappings, the fixed point theorem fails without continuity for every constant.

A $b$-suprametric is a logarithmic transform of a $b$-metric with the same constant, and an interpolative metric is a power of a metric. Consequently, the fixed point theorems of maximum type in these spaces are equivalent to their counterparts in $b$-metric and metric spaces, and linear contractions are Banach contractions for explicit metrics.

The tests of Section~\ref{subsec:tests} make the same analysis available for any axiom given by a triangle function, and Theorem~\ref{thm:char} does the same for any distance that dominates a metric. We suggest applying them before a new axiom is adopted as a setting for fixed point theorems.

\end{document}